\documentclass[11pt]{article}
\usepackage{amsmath,amsfonts,amsthm,mathtools}

 \newtheorem{thm}{Theorem}[section]
 
 \newtheorem{lem}[thm]{Lemma}

 \newtheorem{rem}[thm]{Remark}
\newtheorem{exam}[thm]{Example}

\title{Perturbation bounds for the matrix equation $X+ A^*\widehat{X}^{-1}A=Q$}
\author{Vejdi Hasanov \smallskip\\
\small Faculty of Mathematics and Informatics, \\ 
\small Shumen University, Shumen  9712, Bulgaria,\\
\small e-mail: v.hasanov@shu.bg}
\date{02.03.2018}

\begin{document}
\maketitle

\begin{abstract}
Consider the matrix equation $X+ A^*\widehat{X}^{-1}A=Q$,
where $Q$ is an $n \times n$ Hermitian positive definite matrix, A is an $mn\times n$ matrix, and $\widehat{X}$ is the $m\times m$ block diagonal matrix with $X$ on its diagonal. In this paper, a perturbation bound for the maximal positive definite solution $X_L$
is obtained. Moreover, in case of $\|\widehat{X_L^{-1}}A\|\ge 1$ a modification of the main result is derived. The theoretical results are illustrated by numerical examples.

\bigskip
\noindent Keywords: Nonlinear matrix equation; Positive definite solutions; Maximal solution; Perturbation bounds.
\bigskip \noindent AMS Subject Classification: 15A24; 65H05; 47H14
\end{abstract}

\section{Introduction}

In this paper we study for perturbation bounds the matrix equation
\begin{equation}\label{eq1}
    X+A^*\widehat{X}^{-1}A=Q\,,
\end{equation}
where $Q$ is an $n \times n$ Hermitian positive definite matrix, $A$ is an $mn\times n$ matrix, $\widehat{X}$ is the $m\times m$ block diagonal matrix defined by $\widehat{X}=\mathrm{diag}(X,X,\ldots,X)$, in which $X$ is $n\times n$ matrix, and  $A^*$ is the conjugate transpose of a matrix $A$.

Eq.~(\ref{eq1}) can be write as
\begin{equation}\label{eq1a}
   X +\sum_{i=1}^m A_i^*X^{-1}A_i=Q
\end{equation}
where $A_1, A_2, \ldots, A_m$ are  $n \times n$ matrices, and
 \[A=\left(
  \begin{array}{c}
    A_1 \\
    \vdots \\
    A_m \\
  \end{array}
\right).\]
Moreover, Eq.~(\ref{eq1}) can be reduced to
\begin{equation}\label{eq1i}
  Y+B^*\widehat{Y}^{-1}B=I,
\end{equation}
by multiplying  both hand side of (\ref{eq1}) with the matrix $Q^{-\frac{1}{2}}$, where $I$ is the identity matrix. Thus, Eq.~(\ref{eq1}) is solvable if and only
if Eq.~(\ref{eq1i}) is solvable.

The maximal positive definite solution of Eq.~(\ref{eq1}) with $m=1$ have many applications in ladder networks, control theory, dynamic programming, stochastic filtering, etc., see for instance \cite{AMT,E,Z} and the references therein. Since 1990, the Eq.~(\ref{eq1}) with $m=1$ has been extensively studied, and the research results mainly concentrated on the following: sufficient and necessary conditions for the existence of a positive definite solution \cite{AMT,E,ERR}; numerical methods for computing the positive definite solution \cite{Z,GL,M,IHU}; properties of the positive definite solution \cite{Z,ERR}; and perturbation bounds for the positive definite solution \cite{Xu,SX,HIU,HI,H}.

Eq.~(\ref{eq1i}) is introduced by Long et al. \cite{LHZ} for $m=2$ and by He and Long \cite{HL} for generale case. Later Eqs.~(\ref{eq1}) and (\ref{eq1i}) are investigated by many authors \cite{DLL,PKPA,PPKA,VVSP,HA,HH,HH1,HB}. Bini et al. \cite{BLM} have considered the equation $X+\sum_{i=1}^m C_iX^{-1}D_i = E$ arising in Tree-Like stochastic processes.

Long et al. \cite{LHZ} have given some necessary and sufficient conditions for the existence of a positive definite solution of Eq. (\ref{eq1i}) in case of $m=2$, and proposed basic fixed point iteration and its inversion free variant for finding the largest positive definite solution to that equation. Vaezzadeh et al. \cite{VVSP} have considered inversion free iterative methods for (\ref{eq1}) when $m=2$, also. Hasanov and Ali \cite{HA} improved the results of Vaezzadeh et al. (in \cite{VVSP}) and gave convergence rate of the considered methods. Popchev et al. \cite{PKPA, PPKA} have made a perturbation analysis of (\ref{eq1i}) for $m=2$.

He and Long \cite{HL} have proposed a basic fixed point iteration and its inversion free variant method for finding the maximal positive definite solution to Eq.~(\ref{eq1i}). Hasanov and Hakkaev in \cite{HH} considered the Newton's method for Eq.~(\ref{eq1}) and in \cite{HH1} gave convergence rate of the basic fixed poind iteration and its two inverse free variants, and considered a modification of Newton's method with linear rate of convergence. Duan et al. \cite{DLL} have derived a perturbation bound for the maximal positive definite solution of Eq.~(\ref{eq1i}) based on the matrix differentiation.
Hasanov and Borisova \cite{HB} obtained two perturbed bounds, which do not require the maximal solution to the perturbed or the unperturbed equations.
In addition, many authors have investigated similar or more general nonlinear matrix equations $X-\sum_{i=1}^{m}A_i^*X^{-1}A_i=Q$ \cite{YF,H1}, $X-\sum_{i=1}^{m}A_i^*X^{-\delta_i}A_i=Q$ \cite{DLT,DWL1}, $X+A^*\mathcal{F}(X)A=Q$ \cite{SR, RR}, $X\pm \sum_{i=1}^{m}A_i^*\mathcal{F}(X)A_i=Q$ \cite{RR3}, $X+\sum_{i=1}^{m}A_i^*X^{-q}A_i=Q$ $(0<q\le 1)$ \cite{YWF}, and $A_0+\sum_{i=1}^{k}\sigma_iA_i^*X^{p_i}A_i=0, \ \sigma_i=\pm 1$ \cite{KPPA, PKPAa}.

Motivated by the work in the above papers, we continue to study Eq.~(\ref{eq1}). Here, we derive new perturbation bounds for the maximal solution to Eq.~(\ref{eq1}) by generalization of the results in \cite{HI,H}. Our bounds are much less expensive for computing because they use very simple formulas.

The rest of the paper is organized as follows. In Section~\ref{s2} we give some preliminaries for the perturbation analysis. The main result and some known perturbation bounds are presented in Section~\ref{s3}. Three illustrative examples are provided in Section~\ref{s4}. The paper
closes with concluding remarks in Section~\ref{s5}.

Throughout this paper, we denote by $\mathcal{H}^{n}$ the set
of all $n\times n$ Hermitian matrices. The notation $A>0$ $(A\ge 0)$ means that $A$ is positive definite (semidefinite). If $A-B>0$ (or $A-B\ge 0$) we write $A>B$ (or $A\ge
B$). $I$ (or $I_n$) stands for the identity matrix of order $n$. A Hermitian solution $X_L$ we call maximal one if $X_L\ge X$ for
an arbitrary Hermitian solution $X$. The symbols $\rho(\cdot)$, ${\|\cdot\|}$, ${\|\cdot\|}_F$, and ${\|\cdot\|}_U$ stand the spectral radius, the spectral norm, the Frobenius norm, and any unitary invariant matrix norm, respectively. For $n\times n$ complex matrix $A = (a_{ij})$ and a matrix $B$, $A\otimes B = (a_{ij}B)$ is a Kronecker product. Finally, for a matrix $X$, we denote with $\widehat{Z}$ the $m\times m$ block diagonal matrix with $Z$ on its diagonal, i.e. $\widehat{Z}=I_m\otimes Z$.

\section{Statement of the problem and preliminaries}\label{s2}

It is proved in \cite{HL} that if Eq.~(\ref{eq1i}) has a positive definite solution, then it has a maximal Hermitian solution $X_L$. Moreover, if $\sum_{i=1}^{m}\|B_i\|^2<\frac{1}{4}$, then Eq.~(\ref{eq1i}) with $B=(B_1^T,\ldots, B_m^T)^T$ has maximal positive
definite solution $Y_L$, $\frac{1}{2}I < Y_L \le I$, and it's an unique solution
with these properties. These results are valid for Eq.~(\ref{eq1}) also, i.e.,
if Eq.~(\ref{eq1}) has a positive definite solution, then it has a maximal solution $X_L$. If $\sum_{i=1}^{m}\|Q^{-\frac{1}{2}}A_iQ^{-\frac{1}{2}}\|^2<\frac{1}{4}$, then Eq.~(\ref{eq1}) has maximal positive
definite solution $X_L$, $\frac{1}{2}Q < X_L \le Q$, and it's a unique solution
with these properties.
Moreover, these results have been generalized to equation $X+\sum_{i=1}^{m}A_i^*X^{-q}A_i=Q$ $(0<q\le 1)$ by Yin et al. in \cite{YWF}.

Now, we show that the condition $\sum_{i=1}^{m}\|B_i\|^2<\frac{1}{4}$ for existing of maximal positive
definite solution $Y_L$, for which $\frac{1}{2}I < Y_L \le I$ can be replaced  with $\|B\|<\frac{1}{2}$.

\begin{lem}\label{ll} If $\|B\|<\frac{1}{2}$, then Eq.~(\ref{eq1i}) has a maximal solution $Y_L$ and $\frac{1}{2}I\le Y_L\le I$. Moreover, $\|Y^{-1}\| > 2$ for any other solution $Y$.
\end{lem}

{\bf Proof.} For $C,D\in \mathcal{H}^{n}$, we define a set of matrices as follows
\[[C,D]=\{X \in \mathcal{H}^{n}: C\le X \le D\}.\]
We consider a map $G(Y)=I-B^*\widehat{Y}^{-1}B$. Thus, all the solutions of Eq.~(\ref{eq1i}) are fixed points of $G$. The map $G$ is continuous on $[\frac{1}{2}I, I]$. We prove that $G([\frac{1}{2}I, I])\subset [\frac{1}{2}I, I]$.

Let $Y\in [\frac{1}{2}I, I]$, then
\[I\ge G(Y)=I-B^*\widehat{Y}^{-1}B \ge I - 2 B^*B\ge \frac{1}{2}I.\]

Therefore,  $G(Y)\in[\frac{1}{2}I, I]$ and according to Schauder's fixed point theorem \cite{D} there
exists a matrix $Y_+ \in [\frac{1}{2}I, I]$ such that
$G(Y_+)=Y_+$, i.e., $Y_+$ is a solution of Eq.~(\ref{eq1i}). It is obviously that the maximal solution $Y_L\in [\frac{1}{2}I, I]$. Now, we prove that $Y_L$ is a unique solution in $[\frac{1}{2}I, I]$.

Let $Y_+$ and $Y_L$ are two solution of Eq.~(\ref{eq1i}) in $[\frac{1}{2}I, I]$. We have

\begin{align*}
    0\le\|Y_L-Y_+\|&=\big\|B^*(\widehat{Y_+^{-1}}-\widehat{Y_L^{-1}})B\big\|
   =\big\|B^*\widehat{Y_+^{-1}}(\widehat{Y_L}-\widehat{Y_+})\widehat{Y_L^{-1}}B\big\|\\
   &\le \|B\|^2\big\|\widehat{Y_+^{-1}}\big\|\big\|\widehat{Y_L^{-1}}\big\| \big\|\widehat{Y_L}-\widehat{Y_+}\big\|
    \le 4\|B\|^2\big\|\widehat{Y_L}-\widehat{Y_+}\big\|= 4\|B\|^2\|Y_L-Y_+\|.
\end{align*}

Since $4\|B\|^2<1$, then $Y_L\equiv Y_+$.
\hfill \fbox{}

\begin{rem}\label{r1} By Lemma~\ref{ll} we have, if $\|\widehat{Q^{-\frac{1}{2}}A} Q^{-\frac{1}{2}}\|<\frac{1}{2}$, then Eq.~(\ref{eq1}) has a maximal solution $X_L$ and $\frac{1}{2}Q\le X_L\le Q$. Moreover, $X_L$ is a unique solution in $[\frac{1}{2}Q, Q]$.
\end{rem}

Hasanov and Hakkaev in \cite{HH} have obtained
\begin{equation}\label{roX}
    \rho\Big(\sum_{i=1}^m (X_L^{-1}A_i)^T\otimes (X_L^{-1}A_i)^*\Big) \le 1.
\end{equation}
Moreover, we have (see \cite{H1})
\begin{equation}\label{spN}
    \rho\Big(\sum_{i=1}^m (X_L^{-1}A_i)^T\otimes (X_L^{-1}A_i)^*\Big)
 \le \big\|\widehat{X_L^{-1}}A\big\|^2.
\end{equation}

\begin{lem}\label{l1}\cite{RR2} Let $W\in \mathcal{H}^n$, $C_i$, $i=1,2,\ldots,m$ be $n\times n$ matrices, and $C=(C_1^T,\ldots, C_m^T)^T$. Then
\begin{enumerate}
  \item [(a)] if $\rho(\sum_{i=1}^m C_i^T\otimes C_i^*)<1$, then the equation $X-C^*\widehat{X}C=W$ has a unique solution $P$, and $P\ge 0$ $(P>0)$, when $W\ge 0$ $(W>0)$;
  \item [(b)] if there is some $P>0$ such that $P-C^*\widehat{P}C$ is positive definite, then $\rho(\sum_{i=1}^m C_i^T\otimes C_i^*)<1$.
\end{enumerate}
\end{lem}

\begin{lem}\label{l2} Let $X_+$ be a positive definite solution of Eq.~(\ref{eq1}) with $A=(A_1^T,\ldots, A_m^T)^T$. If
\begin{equation}\label{roX1}
    \rho\Big(\sum_{i=1}^m (X_+^{-1}A_i)^T\otimes (X_+^{-1}A_i)^*\Big) < 1,
\end{equation}
then $X_+\equiv X_L$, i.e., the maximal solution $X_L$ is a unique positive definite solution which satisfy the condition (\ref{roX1}).
\end{lem}

{\bf Proof.} Let $X_+$ be a positive definite solution of Eq.~(\ref{eq1}) which satisfy the condition (\ref{roX1}) and $X_L$ be the maximal solution. Since $X_++A^*\widehat{X_+^{-1}}A=Q$ and $X_L+A^*\widehat{X_L^{-1}}A=Q$, we have
\begin{align*}
    X_+-X_L &= A^*\big(\widehat{X_L^{-1}}-\widehat{X_+^{-1}}\big)A
            = A^*\big(\widehat{X_+^{-1}}+\widehat{X_L^{-1}}-\widehat{X_+^{-1}}\big) \big(\widehat{X_+}-\widehat{X_L}\big)\widehat{X_+^{-1}}A \\
            &= A^*\widehat{X_+^{-1}} \big(\widehat{X_+}-\widehat{X_L}\big)\widehat{X_+^{-1}}A
              + A^*\widehat{X_+^{-1}} \big(\widehat{X_+}-\widehat{X_L}\big)\widehat{X_L^{-1}} \big(\widehat{X_+}-\widehat{X_L}\big)\widehat{X_+^{-1}}A\,.
\end{align*}
Thus,
\[ X_+-X_L - A^*\widehat{X_+^{-1}} \big(\widehat{X_+}-\widehat{X_L}\big)\widehat{X_+^{-1}}A = A^*\widehat{X_+^{-1}} \big(\widehat{X_+}-\widehat{X_L}\big)\widehat{X_L^{-1}} \big(\widehat{X_+}-\widehat{X_L}\big)\widehat{X_+^{-1}}A\,, \]
which implies that $X_+-X_L$ is a solution of the equation $X-C^*\widehat{X}C=W$, where
\begin{align*}
    W &= A^*\widehat{X_+^{-1}} \big(\widehat{X_+}-\widehat{X_L}\big)\widehat{X_L^{-1}} \big(\widehat{X_+}-\widehat{X_L}\big)\widehat{X_+^{-1}}A \ge 0, \\
    C &= \widehat{X_+^{-1}}A, \ \mbox{with} \ C_i=X_+A_i.
\end{align*}
By Lemma~\ref{l2} (a) we have that $X_+-X_L\ge 0$. But, $X_L$ is the maximal solution, i.e. $X_+\le X_L$. Hence, $X_+\equiv X_L$.
\hfill \fbox{}

\begin{rem} We have following hypothesis: the maximal solution $X_L$ is a unique positive definite solution of Eq.~(\ref{eq1}) which satisfy the condition (\ref{roX}).
\end{rem}

\begin{lem}\label{l3} Let $X_+$, be a positive definite solution of Eq.~(\ref{eq1}). If there is a positive definite matrix $P$ such that $\big\|\widehat{PX_+^{-1}}AP^{-1}\big\|<1$, then $X_+$ is a maximal solution, i.e., $X_+\equiv X_L$.
\end{lem}
{\bf Proof.} Let $X_+$, be a positive definite solution of Eq.~(\ref{eq1}) and $P$ is a positive definite matrix $P$ such that $\big\|\widehat{PX_+^{-1}}AP^{-1}\big\|<1$. Then
\begin{align*}
    X_+ +  A^*\widehat{X_+^{-1}}A &= Q, \\
    P^{-1}X_+P^{-1} + P^{-1} A^*\widehat{P}^{-1}\widehat{PX_+^{-1}P}\widehat{P}^{-1}AP^{-1} &= P^{-1}QP^{-1}.
\end{align*}

Therefore, $Y_+=P^{-1}X_+P^{-1}$ is a positive definite solution of the equation
\begin{equation}\label{soseq}
    Y+C^*\widehat{Y^{-1}}C = Q_P
\end{equation}
with $C=\widehat{P}^{-1}AP^{-1}$ and $Q_P=P^{-1}QP^{-1}$.

Since
\[\big\|\widehat{Y_+^{-1}}C\big\|=\big\|\widehat{PX_+^{-1}P}\widehat{P}^{-1}AP^{-1}\big\|
   =\big\|\widehat{PX_+^{-1}}AP^{-1}\big\|<1,\]
by Lemma~\ref{l2} and (\ref{spN}), it follows that $Y_+$ is a maximal solution of Eq.~(\ref{soseq}).
Let $X_L$ be a maximal solution of Eq.~(\ref{eq1}), i.e. $X_L\ge X_+$. Then $Y'=P^{-1}X_LP^{-1}$ is a positive definite solution of Eq.~(\ref{soseq}) and $0 \le Y_+-Y' = P^{-1}(X_+-X_L)P^{-1}$, i.e., $X_+-X_L\ge 0$. Hence, $X_+\equiv X_L$.
\hfill \fbox{}

Consider the perturbed equation
\begin{equation}\label{eq1p}
\tilde X + \tilde A^*\widehat{{\tilde{X}}^{-1}}\tilde A=\tilde Q,
\end{equation}
where $\tilde A = A + \Delta A$, $\tilde Q = Q +\Delta Q$. The matrices $\Delta A$ and $\Delta Q$, ($\Delta Q \in \mathcal{H}^{n}$) are small perturbations in the matrix coefficients $A$ and $Q$ in Eq.~(\ref{eq1}), such that $\tilde Q>0$.

We suppose that Eq.~(\ref{eq1}) has a maximal positive definite solution $X_L$.
The main question is: \textit{how much are the perturbations $\Delta A$ and $\Delta Q$ in the coefficient matrices $A$ and $Q$, respectively such that Eq.~(\ref{eq1p}) has a maximal positive definite solution $\tilde X_L$?} The second question is: \textit{how much is the perturbation $\Delta X_L=\tilde X_L-X_L$, when we have small perturbations $\Delta A$ and $\Delta Q$ in $A$ and $Q$?}

\section{Perturbation bounds}\label{s3}

The questions in the preview section for Eq.~(\ref{eq1}) in case of $m=1$ have been investigated by several authors \cite{Xu,SX,HIU,HI,H}. Hasanov and Ivanov in \cite{HI} have obtained the following result.
\begin{thm}\cite[Theorem 2.1]{HI}\label{t_eq1a}
Let  $A_1, Q, \tilde A_1, \tilde Q$ be
coefficient matrices for equations $X+A_1^*X^{-1}A_1=Q$ and $\tilde X+\tilde A_1^*\tilde X^{-1}\tilde A_1=\tilde Q$. Let
\begin{align*}
  & b_+ = 1-{\|X_L^{-1}A_1\|}^2 + {\|X_L^{-1}\|}\,{\|\Delta Q\|}_U\,, \\
  & c_+ = {\|\Delta Q\|}_U + 2\,{\|X_L^{-1}A_1\|}\,{\|\Delta A_1\|}_U + {\|X_L^{-1}\|}
          \,{\|\Delta A_1\|}_U^2 \,,
\end{align*}
where  $X_L$ is the maximal positive definite solution of the equation  $X+A_1^*X^{-1}A_1=Q$. If
\begin{equation*}\label{condt2}
 {\|X_L^{-1}A_1\|} < 1 \ \ \ \mbox{and} \ \ \ 2\,{\|\Delta{A_1}\|}_U + {\|\Delta{Q}\|}_U
 \le \frac{{(1-{\|X_L^{-1}A_1\|})}^2}{{\|X_L^{-1}\|}} \,,
 \end{equation*}
then $D_+= b_+^2-4c_+{\|X_L^{-1}\|} \ge 0$, the perturbed matrix
equation $\tilde X+\tilde A_1^*\tilde X^{-1}\tilde A_1=\tilde Q$ has the maximal positive definite solution $\tilde X_L$, and
\begin{equation*}\label{cond1}
 {\|\Delta{X_L}\|}_U \ \le \ \frac{b_+-\sqrt{D_+}}{2{\|X_L^{-1}\|}}
   \ =: S_{err}^+ \,.
\end{equation*}
\end{thm}

Moreover, in \cite{HI} has been obtained similar result for equation $X-A^*X^{-1}A=Q$, which was generalized for equation $X-\sum_{i=1}^m A_i^*X^{-1}A_i=Q$ by Yin and Fang \cite{YF}.

Now, we derive new perturbation bounds for the maximal solution to Eq.~(\ref{eq1}) by generalization of Theorem~\ref{t_eq1a} and its modification in \cite{H}. Firstly, we define  $\theta_U(m)={\|\widehat{Z}\|}_U/{\|Z\|}_U$ for an $n\times n$ matrix $Z$ and a unitary invariant norm $\|\cdot\|_U$.
Note that, the values of $\theta_U(m)$ in cases of the spectral norm $\|\cdot\|$ and the Frobenius norm $\|\cdot\|_F$, are $\theta(m)=1$ and $\theta_F(m)=\sqrt{m}$, respectively.
\begin{thm}\label{t_eq1}
Let $A, Q, \tilde A, \tilde Q$ be coefficient matrices for Eqs. (\ref{eq1}) and
(\ref{eq1p}). Let
\begin{align*}
      b&= 1-K_U\big\|\widehat{X_L^{-1}}A\big\| + {\|X_L^{-1}\|}\,{\|\Delta Q\|}_U\,, \\
      c&= {\|\Delta Q\|}_U + 2\,{\big\|\widehat{X_L^{-1}}A\big\|}\,{\|\Delta A\|}_U + {\|X_L^{-1}\|}
          \,{\|\Delta A\|}_U^2 \,,
\end{align*}
where $K_U=\min\big\{\theta_U(m)\big\|\widehat{X_L^{-1}}A\big\|, \big\|\widehat{X_L^{-1}}A\big\|_U\big\}$ and $X_L$ is the maximal positive definite solution of Eq.~(\ref{eq1}). If
\begin{equation}\label{condt2}
 {K_U\big\|\widehat{X_L^{-1}}A\big\|} < 1 \ \ \ \mbox{and} \ \ \ 2\,{\|\Delta{A}\|}_U + {\|\Delta{Q}\|}_U
 < \frac{\Big(1-\sqrt{K_U\big\|\widehat{X_L^{-1}}A\big\|}\, \Big)^2}{{\|X_L^{-1}\|}} \,,
 \end{equation}
then $D= b^2-4 c{\|X_L^{-1}\|} > 0$, the perturbed equation (\ref{eq1p}) has a maximal positive definite solution $\tilde X_L$, and
\begin{equation}\label{cond1}
 {\|\Delta{X_L}\|}_U \ \le \ \frac{b-\sqrt{D}}{2{\|X_L^{-1}\|}}
   \ := \ S_{err} \,.
\end{equation}
\end{thm}
{\bf Proof.} Let  $\tilde X$ be an arbitrary positive definite
solution of Eq.~(\ref{eq1p}). Subtracting
(\ref{eq1}) from (\ref{eq1p}) gives
\begin{equation*}\label{equ_1p}
\Delta X - A^*\widehat{\tilde X^{-1}}\widehat{\Delta X} \widehat{X_L^{-1}} A + A^* \widehat {\tilde X^{-1}}\Delta{A} +(\Delta{A})^*\widehat {\tilde X^{-1}}(A+\Delta A) = \Delta{Q}\,,
\end{equation*}
where  $\Delta{X} = \tilde X - X_L$. Using the equalities
\[\tilde X^{-1} = X_L^{-1} \left(I+\Delta{X} X_L^{-1}\right)^{-1}
        = \left(I+X_L^{-1}\Delta{X}\right)^{-1} X_L^{-1}, \]
we receive
\begin{align}
     \Delta{X} =& \, \Delta{Q} - (\Delta{A})^*
       \left(I_{n^2}+\widehat{X_L^{-1}}\widehat{\Delta{X}}\right)^{-1} \widehat{X_L^{-1}} (A+\Delta{A})   \nonumber\\
    & + A^* \widehat{X_L^{-1}} \left(I_{n^2}+\widehat{\Delta{X}} \widehat{X_L^{-1}}\right)^{-1}
       \left(\widehat{\Delta{X}} \widehat{X_L^{-1}} A - \Delta{A}\right) . \label{equ_2p}
\end{align}

Consider a map $\mu\,:\,\mathcal{H}^n \rightarrow
\mathcal{H}^n$ defined by following way:
\begin{align*}
    \mu(\Delta X) =& \, \Delta{Q}  -  (\Delta{A})^*
       \left(I_{n^2}+\widehat{X_L^{-1}} \widehat{\Delta{X}}\right)^{-1} \widehat{X_L^{-1}} (A+\Delta{A}) \\
    & +  A^* \widehat{X_L^{-1}} \left(I_{n^2}+\widehat{\Delta{X}} \widehat{X_L^{-1}}\right)^{-1}
       \left(\widehat{\Delta{X}} \widehat{X_L^{-1}} A - \Delta{A}\right) .
\end{align*}

Using the inequalities in (\ref{condt2}), we have
\[2 {\|X_L^{-1}\|}{\|\Delta A\|}_U + {\|X_L^{-1}\|}{\|\Delta Q\|}_U
       <  1 + K_U\big\|\widehat{X_L^{-1}}A\big\| - 2 \sqrt{K_U\big\|\widehat{X_L^{-1}}A\big\|}\,, \]
\begin{equation}\label{eqD}
    0 < b < 2 - 2 \Big(\sqrt{K_U\big\|\widehat{X_L^{-1}}A\big\|} + {\|X_L^{-1}\|}{\|\Delta A\|}_U\Big)\,,
\end{equation}
which implies that
\begin{align*}
    D &= b^2 - 4 {\|X_L^{-1}\|} c \\
      &= b^2 - 4 b + 4 - 4\left(K_U{\big\|\widehat{X_L^{-1}}A\big\|} + 2{\|X_L^{-1}\|}{\big\|\widehat{X_L^{-1}}A\big\|}{\|\Delta A\|}_U+\|X_L^{-1}\|^2{\|\Delta A\|}_U^2 \right)\\
      &= (2-b)^2 - 4\Big(\sqrt{K_U\big\|\widehat{X_L^{-1}}A\big\|} + {\|X_L^{-1}\|}{\|\Delta A\|}_U\Big)^2 > 0 \,.
\end{align*}

The square equation
\begin{equation}\label{queqp}
    {\|X_L^{-1}\|}\, S^2 - b \, S + c =0
\end{equation}
has two positive real roots with the smaller one
\[S_{err} = {b- \sqrt{D} \over 2 {\|X_L^{-1}\|}}.\]

We define
\begin{equation*}\label{set1}
    \mathcal{L}_{S_{err}}=\left\{\Delta X \in \mathcal{H}^n: {\|\Delta X\|}_U \le
    S_{err}\right\}.
\end{equation*}

For each  $\Delta X \in \mathcal{L}_{S_{err}}$ we have
\begin{equation}\label{inv-est}
    {\big\|\widehat{X_L^{-1}}\widehat{\Delta X}\big\|}={\|X_L^{-1}\Delta X\|} \le {\|X_L^{-1}\|} {\|\Delta X\|}_U \le {\|X_L^{-1}\|} \, S_{err} < {b\over 2} < 1.
\end{equation}
Thus
$I_{n^2}+\widehat{X_L^{-1}}\widehat{\Delta X}$ is a nonsingular matrix and
\begin{equation}\label{inv-est1}
  \big\|\big(I_{n^2}+\widehat{X_L^{-1}}\widehat{\Delta X}\big)^{-1}\big\|
   \le {1 \over 1 - \|X_L^{-1}\Delta X\|}
   \le {1 \over 1 - {\|X_L^{-1}\|}{\|\Delta X\|}_U}
   \le {1 \over 1 - {\|X_L^{-1}\|} S_{err}} \, .
\end{equation}
According to definition for $\mu(\Delta X)$, for each $\Delta
X \in \mathcal{L}_{S_{err}}$ we obtain
\begin{align*}
     {\|\mu(\Delta X)\|}_U \le&\, {\|\Delta Q\|}_U + {\|\Delta A\|}_U
       \big\| \big(I_{n^2}+\widehat{X_L^{-1}}\widehat{\Delta X}\big)^{-1}
        \widehat{X_L^{-1}}(A+\Delta A) \big\| \\
     &+  \big\|A^* \widehat{X_L^{-1}}
    \big(I_{n^2}+\widehat{\Delta X}\widehat{X_L^{-1}}\big)^{-1}\big\|
    \big\|\big(\widehat{\Delta X} \widehat{X_L^{-1}} A - \Delta A \big)\big\|_U \\
   \le&\, {\|\Delta Q\|}_U + {\|\Delta A\|}_U
       \frac{{\big\|\widehat{X_L^{-1}}A\big\|} + {\|X_L^{-1}\|} {\|\Delta A\|}_U}
       {1 - {\|X_L^{-1}\|} S_{err}}
   + \, {\big\|\widehat{X_L^{-1}}A\big\|} \frac{S_{err}K_U + {\|\Delta A\|}_U}{1 - {\|X_L^{-1}\|} S_{err}} \\
   =&\, \frac{(1-b) S_{err} + c}
       {1 - {\|X_L^{-1}\|} S_{err}} =  S_{err}\,,
\end{align*}
where the last inequality is due to the fact that $S_{err}$ is a
solution of the square equation (\ref{queqp}).

Thus $\mu(\Delta X) \in
\mathcal{L}_{S_{err}}$ for every  $\Delta X \in
\mathcal{L}_{S_{err}}$, which means that $\mu_+(\mathcal{L}_{S_{err}}) \subset \mathcal{L}_{S_{err}}$. Moreover, $\mu_+$ is a continuous  mapping on
$\mathcal{L}_{S_{err}}$. According to Schauder's fixed point theorem \cite{D}  there
exists a $\Delta X_+ \in \mathcal{L}_{S_{err}}$ such that
$\mu(\Delta X_+)= \Delta X_+$. Hence there exists a solution
 $\Delta X_+$ of Eq.~(\ref{equ_2p}) for which
\[{\|\Delta X_+\|}_U \le S_{err}\,.\]

Let
\begin{equation}\label{deftXp}
 \tilde X_+ = X_L + \Delta X_+\,.
\end{equation}
Since  $X_L$ is a solution of  Eq.~(\ref{eq1}) and  $\Delta X_+$ is a
solution of Eq.~(\ref{equ_2p}), then  $\tilde X_+$ is a Hermitian
solution of the perturbed equation  (\ref{eq1p}).

First,  we prove that  $\tilde X_+$ is  a positive definite
solution, and second we prove that  $\tilde X_+ \equiv \tilde
X_L$, i.e, $\tilde X_L \equiv \tilde X_+ = X_L + \Delta X_+$ is
the maximal positive definite solution of Eq.~(\ref{eq1p}).

Since  $X_L$ is a positive definite matrix, then there exists a
positive definite matrix square root of $X_L^{-1}$. From
(\ref{deftXp}) we receive
\[\sqrt{X_L^{-1}}\tilde X_+ \sqrt{X_L^{-1}}
  = I + \sqrt{X_L^{-1}}\Delta X_+\sqrt{X_L^{-1}}\,.\]
Since
 \[\big\|\sqrt{X_L^{-1}}\Delta X_+\sqrt{X_L^{-1}}\big\|
      \le {\|X_L^{-1}\|}{\|\Delta X_+\|}_U < 1\,, \]
then $\sqrt{X_L^{-1}}\tilde X_+ \sqrt{X_L^{-1}} > 0$. Thus, $\tilde
X_+$ is a positive definite solution of Eq.~(\ref{eq1p}). We have to
prove that $\tilde X_+\equiv \tilde X_L$.

Consider $\big\|\widehat{\tilde X_+^{-1}} \tilde A\big\|$. By (\ref{eqD}), (\ref{inv-est}), and (\ref{inv-est1}), we have
\begin{align*}
    \big\|\widehat{\tilde X_+^{-1}} \tilde A\big\| &=
                \big\|\big( I_{n^2} + \widehat{X_L^{-1}}\widehat{\Delta X_+} \big)^{-1}\widehat{X_L^{-1}} (A+\Delta A)\big\| \\
                &\le \frac{\big\|\widehat{X_L^{-1}}A\big\| + \|X_L^{-1}\|\|\Delta A\| }{1-\|X_L^{-1}\|\|\Delta X_+\|}
                \le \frac{\big\|\widehat{X_L^{-1}}A\big\| + \|X_L^{-1}\|{\|\Delta A\|}_U }{1-\|X_L^{-1}\|{\|\Delta X_+\|}_U}\\
                &< \frac{\sqrt{K_U\big\|\widehat{X_L^{-1}}A\big\|} + \|X_L^{-1}\|{\|\Delta A\|}_U }{1- \frac{b}{2}} <1\,.
\end{align*}
Thus, from (\ref{spN}) and Lemma~\ref{l2} (or Lemma~\ref{l3} with $P=I$) it follows that  $\tilde X_+$ is the maximal positive definite solution of Eq.~(\ref{eq1p}), i.e.,  $\tilde X_+\equiv \tilde X_L$
and  $\Delta X_+ \equiv \Delta X_L$.
 \hfill \fbox{}

Note that $\min K_U=\|\widehat{X_L^{-1}}A\|$ and in some cases of Eq.~(\ref{eq1}) the coefficients $A$ and $Q$ have not satisfied the condition $\|\widehat{X_L^{-1}}A\| < 1$.

\begin{exam}\label{ex0} Consider Eq.~(\ref{eq1}) with
\[ A=\left(
      \begin{array}{cc}
        0.1 & 1 \smallskip\\
        1.5 & 10 \smallskip\\
        0.25 & 0.1 \smallskip\\
        0.1 & 1 \\
      \end{array}
    \right)
    \ \ \ \ \mbox{and} \ \ \ \  Q:=X_L+A^*\widehat{X_L^{-1}}A,\]
where
\[ X_L=\left(
            \begin{array}{cc}
              0.5 & -1 \\
              -1 & 50 \\
            \end{array}
          \right)  \]
is the maximal solution.
\end{exam}

For Example~\ref{ex0} we have ${\|\widehat{X_L^{-1}}A\|}=2.5465>1$. Hence, the bound $S_{err}$ in Theorem~\ref{t_eq1} is not applicable. But ${\|\widehat{PX_+^{-1}}AP^{-1}\|}= 0.7316<1$, where $P=\sqrt{Q}$.

\begin{rem} According to Remark~\ref{r1}, from $\big\|\widehat{\sqrt{Q^{-1}}A} \sqrt{Q^{-1}}\big\|<\frac{1}{2}$ it follows
 \[\big\|\widehat{\sqrt{Q}X_+^{-1}}A\sqrt{Q^{-1}}\big\|\le \big\|\sqrt{Q}X_+^{-1}\sqrt{Q}\big\|\big\|\widehat{\sqrt{Q^{-1}}}A\sqrt{Q^{-1}}\big\|< 1.\]
\end{rem}

In case of Example~\ref{ex0}, $\big\|\widehat{\sqrt{Q^{-1}}A} \sqrt{Q^{-1}}\big\|=0.4964$.

Applying the technique developed in \cite{H,H1}, we obtain the following result.
\begin{thm}\label{t_eq1P}
Let  $A, Q, \tilde A, \tilde Q$ be coefficient matrices for Eqs. (\ref{eq1}) and (\ref{eq1p}).
Let
\begin{align*}
  b_p &= 1- K_U^P\big\|\widehat{PX_L^{-1}}AP^{-1}\big\| + \|PX_L^{-1}P\| {\|P^{-1}\Delta QP^{-1}\|}_U\,, \\
  c_p &= {\|P^{-1}\Delta QP^{-1}\|}_U + 2 \big\|\widehat{PX_L^{-1}}AP^{-1}\big\| \big\|\widehat{P^{-1}}\Delta AP^{-1}\big\|_U + \|PX_L^{-1}P\| \big\|\widehat{P^{-1}}\Delta AP^{-1}\big\|_U^2 \,,
\end{align*}
where $K_U^P=\min\big\{\theta_U(m)\big\|\widehat{PX_L^{-1}}AP^{-1}\big\|, \big\|\widehat{PX_L^{-1}}AP^{-1}\big\|_U\big\}$, $X_L$ is the maximal positive definite solution of Eq.~(\ref{eq1}), and $P$ is a positive definite matrix.
If \ $K_U^P\big\|\widehat{PX_L^{-1}}AP^{-1}\big\| < 1$ and
\begin{equation*}\label{condt2a_new}
    2 \big\|\widehat{P^{-1}}\Delta AP^{-1}\big\|_U + {\|P^{-1}\Delta QP^{-1}\|}_U < \frac{\Big(1-\sqrt{K_U^P\big\|\widehat{PX_L^{-1}}AP^{-1}\big\|}\,\Big)^2}{\|PX_L^{-1}P\|} \,,
 \end{equation*}
then $D_p= b_p^2-4c_p\|PX_L^{-1}P\| > 0$ and
\begin{equation}\label{cond1c}
 {\|\Delta{X_L}\|}_U \, \le \, \|P\|^2\frac{b_p-\sqrt{D_p}}{2 \|PX_L^{-1}P\|}
   \, =: \, S_{err}^P \,.
\end{equation}
\end{thm}

\textbf{Proof.} The proof is similar to the proof of Theorem~\ref{t_eq1} by using technique in \cite[Theorem 2.4]{H} and \cite[Theorem 2]{H1}. Moreover, we use Lemma~\ref{l3} for proving that $\tilde X_L = X_L+\Delta X_L$ is a maximal solution of the perturbed equation (\ref{eq1p}). \hfill \fbox{}

Now, we describe some known perturbation bounds.

Xu in \cite{Xu} have obtained an elegant bound in case of $m=1$, which does not require the solution to the perturbed or the unperturbed equations. This bound has been generalized in case of $m>1$ in \cite{HB} and for $Q=I$ in \cite{DLL}.

\begin{thm} \cite[Theorem 3]{HB}\label{HB}
Let
\begin{itemize}
  \item[(i)]  $\displaystyle  \|Q^{-1}\|^2\sum_{i=1}^{m}\|A_i\|^2<\frac{1}{4}$;
  \item[(ii)]  $\displaystyle \|\Delta Q\| \leq \left[\frac{1}{2}-\|Q^{-1}\|\Big( \sum_{i=1}^{m}\|A_i\|^2 \Big)^\frac{1}{2}\right] \|Q^{-1}\|^{-1}$,
  \item[(iii)] $\displaystyle \sum_{i=1}^{m}\|\Delta A_i\|^2  + 2\sum_{i=1}^{m}\|A_i\|\|\Delta A_i\| <\left[\frac{1}{4}-\|\tilde Q^{-1}\|^2\sum_{i=1}^{m}\|A_i\|^2\right]\|\tilde Q^{-1}\|^{-2}$.
\end{itemize}
Then the equations (\ref{eq1}) and (\ref{eq1p}) have maximal solutions $X_L$ and $\tilde{X}_L$,
respectively. Moreover,
\begin{equation*} 
\|\Delta X_L\| \leq \frac{1}{c_1} \left[\|\Delta Q\| +
   2\|\tilde{Q}^{-1}\|\sum_{i=1}^{m}\|\Delta A_i\|(2\|A_i\|+\|\Delta A_i\|)\right]=:est_{hasb17},
\end{equation*}			
where \(c_1 = 1 - 4\|Q^{-1}\|\|\tilde Q^{-1}\|\sum_{i=1}^{m}\|A_i\|^2.\)	
\end{thm}

Theorem~\ref{HB} contains $\|\tilde Q^{-1}\|$. In \cite[Theorem 5]{HB} can be found a perturbation bound which does not depend on the coefficients of the perturbed equation (\ref{eq1p}).

A perturbation bound has been derived for the equation $X+\sum_{i=1}^{m}A_i^*X^{-q}A_i=Q$ $(0<q\le 1)$ by Yin et al. \cite{YWF}. This result rewritten for $q=1$ is as follows

\begin{thm}\cite[Theorem 3.1]{YWF}\label{TYin} Let
\begin{itemize}
  \item[(i)]  $\displaystyle \theta:=\frac{1}{4}-\|Q^{-1}\|^2\sum_{i=1}^m \|A_i\|^2 >0$, \\
  \item[(ii)] $\displaystyle \|\Delta Q\| \le \|Q^{-1}\|^{-1}(1-\sqrt{1-\theta})$, \\
  \item[(iii)] $\displaystyle \sum_{i=1}^m (\|\tilde A_i\|^2 - \|A_i\|^2) <\frac{3}{4} \, \theta \, \|Q^{-1}\|^{-2} $.
\end{itemize}
Then the equations (\ref{eq1}) and  (\ref{eq1p}) have maximal positive definite solutions $X_L$ and $\tilde X_L$, respectively. Moreover,
\begin{equation*}
\|\Delta X_L\|\leq \frac{1}{\xi}\left[\|\Delta Q\| + 2\sum_{i=1}^m \|X_L^{-1}A_i\|\|\Delta A_i\| + \|X_L^{-1}\|\sum_{i=1}^m \|\Delta A_i\|^2\right]=:est_{yinwf14},
\end{equation*}
where
$\xi = 1-{c^2}\sum_{i=1}^m \|\tilde A_i\|^2$ and $c=2\max\{\|Q^{-1}\|, \|\tilde Q^{-1}\|\}$.
\end{thm}

Konstantionov et al. \cite{KPPA} have obtained local and nonlocal perturbation bounds for the equation $A_0+\sum_{i=1}^{k}\sigma_iA_i^*X^{p_i}A_i=0, \ \sigma_i=\pm 1$ by using the technique of Fr'echet derivatives and the method of Lyapunov majorants.
One particular case of this equation is $k=m+1$, $A_0=Q$, $A_{m+1}=I$, $p_{m+1}=1$, $\sigma_{m+1}=\sigma_i=p_i=-1$, $i=1,2,\ldots,m$, i.e. Eq.~(\ref{eq1}).

Now, we formulate the results from \cite{KPPA} in this particular case.
We use some notations. Let
\begin{align}
 & \delta = \big(\|\Delta Q\|_F, \|\Delta A_1\|_F, \cdots, \|\Delta A_m\|_F\big)^T, \label{de}\\
 &L=I_{n^2}-\sum_{i=1}^m (X_+^{-1}A_i)^{T} \otimes {(X_+^{-1}A_i)}^*, \nonumber\\
 &\Pi=\sum_{j=1}^n \sum_{k=1}^n e_j e_k^T\otimes e_k e_j^T, \nonumber
\end{align}
where $e_j$ denotes the $j$th column of $I_n$.

Let
\begin{align*}
   & W_Q=L^{-1}=W_{Q0}+\mathrm{i}W_{Q1}, \nonumber\\
   & W_{A_i} = - L^{-1}(I_n\otimes (X_+^{p_i}A_i)^*) = W_{A_i0}+\mathrm{i}W_{A_i1}, \\
   & W_{\bar A_i}= - L^{-1}((X_+^{p_i}A_i)^T \otimes I_n)\Pi = W_{\bar A_i0}+\mathrm{i}W_{\bar A_i1}, \nonumber\\
   & M_{A_i}=\left(
          \begin{array}{cc}
            W_{A_i0}+ W_{\bar A_i0} & W_{\bar A_i1}-W_{A_i1} \\
            W_{\bar A_i1}+W_{A_i1}  & W_{A_i0}- W_{\bar A_i0} \\
          \end{array}
        \right), \nonumber\\
   & W_{Q}^{\mathcal R}=\left(
          \begin{array}{cc}
            W_{Q0} & -W_{Q1} \\
            W_{Q1}  & W_{Q0} \\
          \end{array}
        \right), \ \  k_{Q}=\|W_Q\|, \nonumber\\
   & k_{A_i}=\|M_{A_i}\|, \ \ \ i=1,2,\ldots, m, \nonumber\\
   & \Gamma =\big(\Gamma_1,\Gamma_2,\ldots,\Gamma_{m+2}\big)= \big(W_{Q}^{\mathcal R},M_{A_1},\ldots,M_{A_{m+1}}\big).\nonumber
\end{align*}

Konstantinov et al. \cite{KPPA} have obtained the local perturbation bounds:
\begin{align*}
   & est_1(\delta) = k_{Q}\|\Delta Q\|_F + \sum_{i=1}^m k_{A_i}\|\Delta A_i\|_F, \\
   & est_2(\delta) = \|\Gamma\|\|\delta\|, \ \ \
   est_3(\delta) = \sqrt{\delta^TR\delta}, \\
   & est(\delta) = \min\{est_2(\delta), est_3(\delta)\},
\end{align*}
where $R$ is an $(m+1)\times (m+1)$ real symmetric matrix with non-negative entries $r_{ij}=\|\Gamma_i^T\Gamma_j\|$, $i,j=1,2,\ldots,m+1$.

We note that, in case of real matrix coefficients in Eq.~(\ref{eq1}), the above formulas are more simple (see \cite{KPPA}).

Let
\begin{align}
    & a_0(\delta) = est(\delta) + \|L^{-1}\|\|X_+^{-1}\|\sum_{i=1}^m\|\Delta A_i\|_F^2, \label{a0}\\
    & a_1(\delta) = \|L^{-1}\| \|X_+^{-1}\|^2\sum_{i=1}^m (2\|A_i\|+\|\Delta A_i\|_F)\|\Delta A_i\|_F, \label{a1}\\
    & a_2(\delta) = \|L^{-1}\|  \|X_+^{-1}\|^3\sum_{i=1}^m (\|A_i\|+\|\Delta A_i\|_F)^2, \label{a2}\\
    & \Omega=\left\{\delta \ \mathrm{from} \ (\ref{de}) :  a_1(\delta) + 2\sqrt{a_0(\delta)a_2(\delta)}\le 1 \right\}. \label{Om}
\end{align}
The following non-local perturbation bound was obtained in \cite{KPPA}.

\begin{thm}{\rm (\cite[Theorem 5.1]{KPPA})} Let $\delta\in \Omega$, where $\Omega$ is given in (\ref{Om}). Then the non-local perturbation bound
\begin{equation}\label{pKPPA}
    \|\Delta X_+\|_F \le \frac{2a_0(\delta)}{1-a_1(\delta)+\sqrt{(1-a_1(\delta))^2-4a_0(\delta)a_2(\delta)}}=:est_{konppa11}
\end{equation}
 is valid for Eq.~(\ref{eq1}), where $a_i(\delta)$, $i=0,1,2$ are determined by (\ref{a0})-(\ref{a2}).
\end{thm}

\section{Numerical experiments}\label{s4}

We experiment with our bounds and the corresponding
perturbation estimates proposed by Hasanov and Borisova \cite{HB}, Yin et al. \cite{YWF} for the equation $X+\sum_{i=1}^{m}A_i^*X^{-q}A_i=Q$ $(0<q\le 1)$ and Konstantionov et al. \cite{KPPA} for the equation $A_0+\sum_{i=1}^{k}\sigma_iA_i^*X^{p_i}A_i=0, \ \sigma_i=\pm 1$.
Denote the ratio of the perturbation bounds to the estimated value as follows:
\begin{align*}
    & hasb17=\frac{est_{hasb17}}{\|\Delta X_L\|}, \hspace{6mm}
    yinwf14=\frac{est_{yinwf14}}{\|\Delta X_L\|}, \hspace{6mm}
    konppa11=\frac{est_{konppa11}}{\|\Delta X_L\|_F},   \\
    & has=\frac{S_{err}}{\|\Delta X_L\|}, \hspace{6mm} has_F=\frac{S_{err}}{\|\Delta X_L\|_F}, \hspace{6mm}
    hasP=\frac{S_{err}^{P}}{\|\Delta X_L\|}, \hspace{6mm} has_FP=\frac{S_{err}^{P}}{\|\Delta X_L\|_F},
\end{align*}
where for $has$ and $hasP$ the perturbation bounds $S_{err}$ and $S_{err}^{P}$ are computed by using the spectral norm, and for $has_F$ and $has_FP$, $S_{err}$ and $S_{err}^{P}$ are computed by using the Frobenius norm. Moreover, we compute $S_{err}^{P}$ for different $P$: $P=\sqrt{Q}$, and $P=P_1:=\sqrt{Q}+4\sqrt[4]{Q}$.

\begin{exam}\label{ex1} Consider Eq.~(\ref{eq1}) with matrices
\[A=\left(
      \begin{array}{c}
        A_1 \\
        A_2 \\
      \end{array}
    \right), \ \ \ Q:=X_L+A^*\widehat{X_L^{-1}}A,
\]
where $ X_L=\mathrm{diag}(0.725, 2, 3, 2, 1)$ is the maximal solution, and
\[ A_1=\frac{2\sqrt{3}}{45}\left(
      \begin{array}{ccccc}
       1 & 0 & 0 & 0 & 1 \smallskip\\
      -1 & 1 & 0 & 0 & 1 \smallskip\\
      -1 &-1 & 1 & 0 & 1 \smallskip\\
      -1 &-1 &-1 & 1 & 1 \smallskip\\
      -1 &-1 &-1 &-1 & 1 \\
      \end{array}
    \right),  \ \ \
  A_2=\frac{1}{15}\left(
      \begin{array}{ccccc}
       2 & 0 & 0 & 0 & 0 \smallskip\\
       1 & 2 & 0 & 0 & 0 \smallskip\\
       1 & 1 & 2 & 0 & 0 \smallskip\\
       0 & 1 & 1 & 2 & 0 \smallskip\\
       0 & 0 & 1 & 1 & 2 \\
      \end{array}
    \right),
     \]
Assume that the
perturbations on $A$ and $Q$ are
\[\Delta A=10^{-2j}\left(
\begin{array}{c}
  2I-0.5E \\
  \|C\|^{-1}C
\end{array}\right), \ \
 \Delta Q=\tilde X_L+ \tilde A^*\widehat{\tilde X_L^{-1}}\tilde A - Q\,,
\]
where $\tilde X_L=X_L+\frac{10^{-2j}}{2}(I-E)$, $\tilde A=A+\Delta A$, $C$ is a random
matrix, which is generated by Matlab's function \textsf{randn}, and $E$ being the $5 \times 5$ matrix with all entries equal to one.
\end{exam}

\begin{table}[h]
 \caption{Results for Example~\ref{ex1}} \label{t1}
\begin{tabular}{crrrr}
  \hline
                   & $j=2$         & $j=3$         & $j=4$         & $j=5$  \\ \hline
 $\|\Delta X_L\|$  & $2.0000e-04$  & $2.0000e-06$  & $2.0000e-08$  & $2.0000e-10$ \\
 $\|\Delta X_L\|_F$& $2.2361e-04$  & $2.2361e-06$  & $2.2361e-08$  & $2.2361e-10$ \\ \hline
 $konppa11$        & $    2.1430$  & $    2.2136$  & $    2.1056$  & $    2.0771$ \\
 $yinwf14$         & $    6.9245$  & $    7.2893$  & $    6.8848$  & $    6.6668$ \\
 $hasb17$          & $   12.7505$  & $   13.1278$  & $   12.7234$  & $   12.5055$ \\
 $has$             & $    1.6571$  & $    1.7470$  & $    1.6308$  & $    1.5983$ \\
 $has_F$           & $    2.3276$  & $    2.4001$  & $    2.3079$  & $    2.2717$ \\
  \hline
\end{tabular}
\end{table}

The ratio of the perturbation bounds and the estimated value for $j=2,3,4,5$ are listed in Table~\ref{t1}. Among the bounds considered in this example the bound $S_{err}$ by using spectral norm, followed by $est_{konppa11}$ and $S_{err}$ by using Frobenius norm, gives the
sharpest estimates. The bound $est_{hasb17}$ is too conservative, but it does not require the solution to the perturbed or the unperturbed equations.

\begin{exam}\label{ex2} Consider Eq.~(\ref{eq1}) with matrices
\[A=\left(
      \begin{array}{c}
        A_1 \\
        A_2 \\
        A_3 \\
      \end{array}
    \right), \ \ \ Q:=X_L+A^*\widehat{X_L^{-1}}A,
\]
where $A_1=\frac{1+\mathrm{i}}{25}A_0$, $A_2=\frac{1+\mathrm{i}}{25}A_0^T$, $A_3=\frac{1}{70}A_0^TA_0$, $\mathrm{i}=\sqrt{-1}$, $X_L=E+1.5I$ is the maximal solution, and
\[ A_0=\left(
      \begin{array}{ccccc}
       1 & 0 & 0 & 0 & 1 \smallskip\\
      -1 & 1 & 0 & 0 & 1 \smallskip\\
      -1 &-1 & 1 & 0 & 1 \smallskip\\
      -1 &-1 &-1 & 1 & 1 \smallskip\\
      -1 &-1 &-1 &-1 & 1 \\
      \end{array}
    \right).
     \]
Assume that the perturbations on $A$ and $Q$ are
\[\Delta A= 10^{-2j}\left(
\begin{array}{c}
  15\mathrm{i}C_0 \\
  25\mathrm{i}C_0^T \\
  C_0^T+C_0 \end{array}\right), \ \
 \Delta Q=\tilde X_L+ \tilde A^*\widehat{\tilde X_L^{-1}}\tilde A - Q\,,
\]
where $C_0=\|C\|^{-1}C$, $\tilde X_L=X_L-10^{-2j}(I+0.25E)$, $\tilde A=A+\Delta A$, and $C$ is a random matrix.
\end{exam}

The ratio of the perturbation bounds and the estimated value for $j=2,3,4,5$ are listed in Table~\ref{t2}. The results for Example~\ref{ex2} are identical with these of  Example~\ref{ex1}.

\begin{table}[h]
 \caption{Results for Example~\ref{ex2}} \label{t2}
\begin{tabular}{crrrr}
  \hline
                   & $j=2$         & $j=3$         & $j=4$         & $j=5$  \\ \hline
 $\|\Delta X_L\|$  & $2.2500e-04$  & $2.2500e-06$  & $2.2500e-08$  & $2.2500e-10$ \\
 $\|\Delta X_L\|_F$& $3.0104e-04$  & $3.0104e-06$  & $3.0104e-08$  & $3.0104e-10$ \\ \hline
 $konppa11$        & $    4.7381$  & $    6.0984$  & $    5.2278$  & $    5.2508$ \\
 $yinwf14$         & $    6.0738$  & $    6.3788$  & $    6.1441$  & $    6.6765$ \\
 $hasb17$          & $   11.8124$  & $   12.1426$  & $   11.9011$  & $   12.4454$ \\
 $has$             & $    4.5401$  & $    5.0005$  & $    4.6716$  & $    5.2315$ \\
 $has_F$           & $    4.6747$  & $    6.0147$  & $    5.1367$  & $    5.1706$ \\
  \hline
\end{tabular}
\end{table}

\begin{exam}\label{ex3} Consider Eq.~(\ref{eq1}) with matrices $A$ and $Q$ from Example~\ref{ex0}, i.e., \[ A=\left(
      \begin{array}{cc}
        0.1 & 1 \smallskip\\
        1.5 & 10 \smallskip\\
        0.25 & 0.1 \smallskip\\
        0.1 & 1 \\
      \end{array}
    \right)
    \ \ \ \ \mbox{and} \ \ \ \  Q:=X_L+A^*\widehat{X_L^{-1}}A,\]
where
\[ X_L=\left(
            \begin{array}{cc}
              0.5 & -1 \\
              -1 & 50 \\
            \end{array}
          \right)  \]
is the maximal solution. Assume that the
perturbations on $A$ and $Q$ are
\[\Delta A=\left(\begin{array}{cc} 1 & 2 \\ 3 & 4 \end{array}\right)
\times 10^{-5}, \ \ \ \
 \Delta Q=\left(\begin{array}{cc} 1 & 5 \\5 & 4 \end{array}\right)\times 10^{-10}\,.\]
\end{exam}

\begin{table}[h]
 \caption{Results for Example~\ref{ex3}} \label{t3}
\begin{tabular}{crrrr}
  \hline
                   & $j=5$         & $j=6$         & $j=7$         & $j=8$  \\ \hline
 $\|\Delta X_L\|$  & $9.6745e-05$  & $9.6743e-06$  & $9.6743e-07$  & $9.6743e-08$ \\
 $\|\Delta X_L\|_F$& $1.0250e-04$  & $1.0249e-05$  & $1.0249e-06$  & $1.0249e-07$ \\ \hline
 $konppa11$        & $        * $  & $        * $  & $        * $  & $   34.5061$ \\
 $yinwf14$         & $        * $  & $        * $  & $        * $  & $        * $ \\
 $hasb17$          & $        * $  & $        * $  & $        * $  & $        * $ \\
 $has$             & $        * $  & $        * $  & $        * $  & $        * $ \\
 $has\sqrt{Q}$     & $  130.8342$  & $  130.7676$  & $  130.7609$  & $  130.7603$ \\
 $hasP_1$            & $   42.9846$  & $   42.9368$  & $   42.9320$  & $   42.9316$ \\
 $has_F$           & $        * $  & $        * $  & $        * $  & $        * $ \\
 $has_F\sqrt{Q}$   & $  133.5516$  & $  133.4700$  & $  133.4619$  & $  133.4611$ \\
 $has_FP_1$        & $   54.4271$  & $   54.3154$  & $   54.3043$  & $   54.3032$ \\
 \hline
\end{tabular}
\end{table}

We recall  that for Example~\ref{ex3} (see Example~\ref{ex0}) we have ${\|\widehat{X_L^{-1}}A\|}=2.5465>1$. Hence, the bound $S_{err}$ in Theorem~\ref{t_eq1} is not applicable. But ${\big\|\widehat{\sqrt{Q}X_+^{-1}}A\sqrt{Q^{-1}}\big\|}= 0.7316$. Moreover,
${\big\|\widehat{\sqrt{Q}X_+^{-1}}A\sqrt{Q^{-1}}\big\|_F}=  0.7742$,  ${\|\widehat{P_1X_+^{-1}}AP_1^{-1}\|}= 0.8089<1$, and ${\|\widehat{P_1X_+^{-1}}AP_1^{-1}\|_F}= 0.9056$, where $P_1=\sqrt{Q}+4\sqrt[4]{Q}$. The ratio of the perturbation bounds and the estimated value for $j=5,6,7,8$ are listed in Table~\ref{t3}. The cases when the conditions of existence of a bound are violated are denoted by an asterisk.

\section{Concluding remarks}\label{s5}

Analyzing the behaviour of the perturbation bounds considered in the paper, we can
point out as most effective the bounds $S_{err}$ and $est_{konppa11}$. When $\|\widehat{X_L^{-1}}A\|>1$ we use the bound $S_{err}^P$ with appropriate matrix $P$. The optimal choosing of matrix $P$ is an open problem. The perturbation bounds  $S_{err}$ or $S_{err}^P$, derived in this paper can be easily computed using any unitary invariant norm $\|\cdot\|_U$, while the bound $est_{konppa11}$ depends on many parameters, which is very difficult for computing in generally. The bound $est_{hasb17}$ is an a priori estimate, since for its calculation it is not necessary to know the solutions $X_L$ and ${\tilde X}_L$ of the unperturbed and the perturbed equation, respectively.

\section*{Acknowledgment}
\thanks{This research work was supported by the Shumen University under Grant No RD-08-145/2018.}

\end{document}